\documentclass[a4paper,11pt]{article}
\usepackage[top=3cm,bottom=3cm,left=2.5cm,right=2.5cm]{geometry}
\usepackage{amsfonts,epsf,amsmath,amssymb,here}
\usepackage{latexsym,multirow,rotating}
\usepackage{pgf,tikz,subfigure}

\newtheorem{theorem}{\bf Theorem}[section]

\newtheorem{proposition}[theorem]{\bf Proposition}

\newcommand{\proof}{\noindent{\bf Proof.\ }}
\newcommand{\qed}{\hfill $\square$ \bigskip}

\begin{document}

\baselineskip=0.30in

\vspace*{40mm}

\begin{center}
{\LARGE \bf\sc $M$-Polynomial and Degree-Based Topological Indices}
\bigskip \bigskip

{\large \sc Emeric Deutsch}\\
\smallskip
{\em Polytechnic Institute of New York University, United States} \\
e-mail: {\tt emericdeutsch@msn.com}

\bigskip

{\large \sc Sandi Klav\v zar} \\
\smallskip
{\em Faculty of Mathematics and Physics, University of Ljubljana, Slovenia} \\
{\em Faculty of Natural Sciences and Mathematics, University of Maribor, Slovenia} \\
{\em Institute of Mathematics, Physics and Mechanics, Ljubljana, Slovenia} \\
e-mail: {\tt sandi.klavzar@fmf.uni-lj.si}

\bigskip\medskip
(Received *** July, 2014)
\end{center}

\noindent
{\bf Abstract}

\vspace{3mm}\noindent

Let $G$ be a graph and let $m_{ij}(G)$, $i,j\ge 1$, be the number of edges $uv$ of $G$ such that $\{d_v(G), d_u(G)\} = \{i,j\}$. The {\em $M$-polynomial} of $G$ is introduced with $\displaystyle{M(G;x,y) = \sum_{i\le j} m_{ij}(G)x^iy^j}$. It is shown that degree-based topological indices can be routinely computed from the polynomial, thus reducing the problem of their determination in each particular case to the single problem of determining the $M$-polynomial. The new approach is also illustrated with  examples. 

\vspace{5mm}

\baselineskip=0.30in

% \noindent {\bf Key words}: 
% 
% \bigskip\noindent
% {\bf AMS subject classification (2010)}: 05C07, 05C31, 92E10

%%%%%%%%%%%%%%%%%%%%%%%%%%%%%%%%%%%%%%%%%%%%%%%%%%%%%%%%%%%%%%%%%%%%%%%%%%%%%
\section{Introduction}
%%%%%%%%%%%%%%%%%%%%%%%%%%%%%%%%%%%%%%%%%%%%%%%%%%%%%%%%%%%%%%%%%%%%%%%%%%%%%

Numerous graph polynomials were introduced in the literature, several of them turned out to be applicable in mathematical chemistry. For instance, the Hosoya polynomial~\cite{hosoya-1988}, see also~\cite{deutsch-2013,eliasi-2013,lin-2013}, is the key polynomial in the area of distance-based topological indices. In particular, the Wiener index can be computed as the first derivative of the Hosoya polynomial, evaluated at 1, the hyper-Wiener index~\cite{cash-2002} and the Tratch-Stankevich-Zefirov index can be obtained similarly~\cite{bruckler-2011,gutman-2012}. Additional chemically relevant polynomials are the matching polynomial~\cite{gutman-1977,farell-1979}, the Zhang-Zhang polynomial (also known as the Clar covering polynomial)~\cite{zhang-1996,chou-2014,zhang-2013}, 
the Schultz polynomial~\cite{hossaini-2013}, and the Tutte polynomial~\cite{doslic-2013}, to name just a few of them. In this paper we introduce a polynomial called the $M$-polynomial, and show that its role for    degree-based invariants is parallel to the role of the Hosoya polynomial for distance-based invariants. 

In chemical graph theory (too) many topological indices were introduced, in particular (too) many degree-based topological indices. This fact is emphasized in the recent survey~\cite{gutman-2013a} which contains a uniform approach to the degree-based indices and a report on a comparative test from~\cite{gutman-2013b} how these indices are correlated with physico-chemical parameters of octane isomers. The test indicates that quite many of these indices are inadequate for any structure-property correlation. But it could be that they are still applicable in a combination with other indices. 

In the literature one finds many papers that, for a given family of graphs, determine a closed formula for a given (degree-based) topological index. To overcome this particular approach in the area of degree-based indices, in this paper we introduce the $M$-polynomial and demonstrate that in numerous cases a degree-based topological index can be expressed as a certain derivative or integral (or both) of the corresponding $M$-polynomial. This in particular implies that knowing the $M$-polynomial of a given family of graphs, a closed formula for any such index can be obtained routinely. In the remaining cases when the function defining a given topological index is such that it does not allow it to be (routinely) determined  from the $M$-polynomial, one can use equality~\eqref{eq:degree-index-2} from the following section and try to get a closed formula from it. But in any case, knowing the $M$-polynomial is sufficient, hence the possible future research in the area should concentrate on determining the $M$-polynomial of a relevant family of graphs instead of computing the corresponding indices one by one. Moreover, it is our hope that a closer look to the properties of the $M$-polynomial will bring new general insights.  

In the next section the $M$-polynomial is introduced and shown how degree-based indices can be computed from it. In the last section it is discussed how the $M$-polynomial can be computed for (chemical) graphs and typical examples for the use of the new polynomial are presented. 

%%%%%%%%%%%%%%%%%%%%%%%%%%%%%%%%%%%%%%%%%%%%%%%%%%%%%%%%%%%%%%%%%%%%%%%%%%%%%
\section{The $M$-polynomial}
%%%%%%%%%%%%%%%%%%%%%%%%%%%%%%%%%%%%%%%%%%%%%%%%%%%%%%%%%%%%%%%%%%%%%%%%%%%%%

If $G=(V,E)$ is a graph and $v\in V$, then $d_v(G)$ (or $d_v$ for short if $G$ is clear from the context) denotes the degree of $v$. Let $G$ be a graph and let $m_{ij}(G)$, $i,j\ge 1$, be the number of edges $e=uv$ of $G$ such that $\{d_v(G), d_u(G)\} = \{i,j\}$. As far as we know, the quantities $m_{ij}$ were first introduced and applied in~\cite{gutman-2002}. We now introduce the {\em $M$-polynomial} of $G$ as
$$M(G;x,y) = \sum_{i\le j} m_{ij}(G)x^iy^j\,.$$
For a graph $G=(V,E)$, a {\em degree-based topological index} is a graph invariant of the form 
\begin{equation}
\label{eq:degree-index}
I(G) = \sum_{e=uv\in E}f(d_u, d_v)\,,
\end{equation}
where $f=f(x,y)$ is a function appropriately selected for possible chemical applications~\cite{gutman-2013a}. For instance, the generalized Randi\'c index $R_\alpha(G)$, $\alpha\ne 0$, is defined with~\eqref{eq:degree-index} by setting $f(x,y) = (xy)^\alpha$~\cite{bollobas-1998}. Collecting edges with the same set of end-degrees we can rewrite~\eqref{eq:degree-index} as
\begin{equation}
\label{eq:degree-index-2}
I(G) = \sum_{i\le j}m_{ij}(G) f(i, j)\,,
\end{equation}
cf.~\cite{deng-2014,deng-2011,gutman-2013a,vu-2009}.  

Using the operators $D_x$ and $D_y$ defined on differentiable functions in two variables by  
$$D_x(f(x,y)) = x \frac{\partial f(x,y)}{\partial x}, \quad 
D_y(f(x,y)) = y \frac{\partial f(x,y)}{\partial y}\,,$$ 
we have: 

\begin{theorem}
\label{thm:main}
Let $G=(V,E)$ be a graph and let $\displaystyle{I(G) = \sum_{e=uv\in E}f(d_u, d_v)}$, where $f(x,y)$ is a polynomial in $x$ and $y$. Then 
$$I(G) = f(D_x,D_y)(M(G;x,y))\big|_{x=y=1}\,.$$
\end{theorem}

\proof
Let $r\ge 0$ and $s\ge 0$ be fixed. Then  
\begin{eqnarray*}
D_x^r D_y^s (M(G;x,y)) & = & 
 D_x^r D_y^s \left(\sum_{i\le j}m_{ij}(G) x^i y^j\right) 
 = D_x^r \left(\sum_{i\le j}m_{ij}(G) j^s x^i y^j\right) \\
 & = & \sum_{i\le j}m_{ij}(G) i^r j^s x^i y^j\,. 
\end{eqnarray*}
Suppose that $f(x,y) = \sum_{r,s}\alpha_{r,s} x^r y^s$. Then 
$f(D_x,D_y)(M(G;x,y))$ is equal to  
$$\sum_{r,s} \alpha_{r,s}D_x^r D_y^s (M(G;x,y)) 
=  \sum_{r,s} \alpha_{r,s} \sum_{i\le j}m_{ij}(G) i^r j^s x^i y^j 
= \sum_{i\le j} m_{ij}(G) \sum_{r,s} \alpha_{r,s}  i^r j^s x^i y^j\,.$$
It follows that 
$$f(D_x,D_y)(M(G;x,y))\big|_{x=y=1} = 
\sum_{i\le j} m_{ij}(G) \sum_{r,s} \alpha_{r,s}  i^r j^s = 
\sum_{i\le j} m_{ij}(G) f(i,j)\,.$$
The result follows by~\eqref{eq:degree-index-2}. 
\qed

In addition to the operators $D_x$ and $D_y$ we will also make use of the operators 
$$S_x(f(x,y)) = \int_{0}^x \frac{f(t,y)}{t} dt, \quad
S_y(f(x,y)) = \int_{0}^y \frac{f(x,t)}{t} dt\,.$$ 
As in the proof of Theorem~\ref{thm:main} we infer that the operators $S_x^r S_y^s$, $S_x^r D_y^s$, and $D_x^r S_y^s$ applied to the term $x^iy^j$ of $M(G;x,y)$, and evaluated at $x=y=1$, return $i^{-r} j^{-s}$, $i^{-r} j^s$, and $i^r j^{-s}$, respectively. Hence, proceeding along the same lines as in the proof of Theorem~\ref{thm:main}, we can state the following extension of Theorem~\ref{thm:main}. 

\begin{theorem}
\label{thm:main2}
Let $G=(V,E)$ be a graph and let $\displaystyle{I(G) = \sum_{e=uv\in E}f(d_u, d_v)}$, where 
$$\displaystyle{f(x,y) = \sum_{i,j\in \mathbb{Z}}\alpha_{ij}x^iy^j}\,.$$ 
Then $I(G)$ can be obtained from $M(G;x,y)$ using the operators $D_x$, $D_y$, $S_x$, and $S_y$. 
\end{theorem}

In the statement of Theorem~\ref{thm:main2}, $I(G)$ appears implicitly (contrary to the explicit formula of Theorem~\ref{thm:main}). Hence in Table~\ref{table1} explicit expressions for some frequent degree-based topological indices are listed. The general case for any function of the form $f(x,y) = \sum_{i,j\in \mathbb{Z}}\alpha_{ij}x^iy^j$ should then be clear. Note that for the symmetric division index we can apply Theorem~\ref{thm:main2} because its defining function $f(x,y)$ can equivalently be written as $x/y + y/x$. 

\begin{table}[ht!]
\begin{center}
\begin{tabular}{|c|c|c|}\hline
$\phantom{\Big|}$topological index$\phantom{\Big|}$ & $f(x,y)$ & derivation from $M(G;x,y)$ \\
\hline\hline
$\phantom{\Big|}$first Zagreb$\phantom{\Big|}$ & 
$x+y$ & $(D_x + D_y)(M(G;x,y))\big|_{x=y=1}$ \\
\hline
$\phantom{\Big|}$second Zagreb$\phantom{\Big|}$ & 
$xy$ & $(D_xD_y)(M(G;x,y))\big|_{x=y=1}$ \\
\hline
$\phantom{\Big|}$second modified Zagreb$\phantom{\Big|}$ & 
$\frac{1}{xy}$ & $(S_xS_y)(M(G;x,y))\big|_{x=y=1}$ \\
\hline
$\phantom{\Big|}$general Randi\'c ($\alpha\in {\mathbb{N}}$)$\phantom{\Big|}$ & 
$(xy)^\alpha$ & $(D_x^\alpha D_y^\alpha)(M(G;x,y))\big|_{x=y=1}$ \\
\hline
$\phantom{\Big|}$general Randi\'c ($\alpha\in {\mathbb{N}}$)$\phantom{\Big|}$ & 
$\frac{1}{(xy)^\alpha}$ & $(S_x^\alpha S_y^\alpha)(M(G;x,y))\big|_{x=y=1}$ \\
\hline
$\phantom{\Big|}$symmetric division index$\phantom{\Big|}$ & 
$\frac{x^2 + y^2}{xy}$ & $(D_xS_y + D_yS_x)(M(G;x,y))\big|_{x=y=1}$ \\
\hline
\end{tabular}
\end{center}
\caption{Some standard degree based topological indices and the formulas how to compute them from the $M$-polynomial}
\label{table1}
\end{table}

In order to handle additional indices that are not covered by Theorem~\ref{thm:main2}, we introduce two additional operators as follows: 
$$J(f(x,y)) = f(x,x), \quad 
Q_\alpha(f(x,y)) = x^\alpha f(x,y), \alpha\ne 0\,.$$ 
With these two operators we have the following result whose proof again proceeds along the same lines as the proof of Theorem~\ref{thm:main}: 

\begin{proposition}
\label{prp:xxx}
Let $G=(V,E)$ be a graph and let $\displaystyle{I(G) = \sum_{e=uv\in E}f(d_u, d_v)}$, where 
$f(x,y) = \frac{ x^r y^s}{(x+y+\alpha)^k}$, where $r,s\ge 0$, $,t\ge 1$, and $\alpha\in {\mathbb{Z}}$. Then  
$$I(G) = S_x^k\,Q_\alpha\,J\,D_x^r\,D_y^s (M(G;x,y))\big|_{x=1}\,.$$ 
\end{proposition}

In Table~\ref{table2} three special cases of Proposition~\ref{prp:xxx} from the literature are collected.  

\begin{table}[ht!]
\begin{center}
\begin{tabular}{|c|c|c|}\hline
$\phantom{\Big|}$topological index$\phantom{\Big|}$ & $f(x,y)$ & derivation from $M(G;x,y)$ \\
\hline\hline
$\phantom{\Big|}$harmonic$\phantom{\Big|}$ & 
$\frac{2}{x+y}$ & $2\,S_x\,J\,(M(G;x,y))\big|_{x=1}$ \\
\hline
$\phantom{\Big|}$inverse sum$\phantom{\Big|}$ & 
$\frac{xy}{x+y}$ & $S_x\,J\,D_x\,D_y\,(M(G;x,y))\big|_{x=1}$ \\
\hline
$\phantom{\Big|}$augmented Zagreb$\phantom{\Big|^{X^X}}$ & $\left(\frac{xy}{x+y-2}\right)^3$ $\phantom{\Big|}$ & $S_x^3\,Q_{-2}\,J\,D_x^3\,D_y^3\,(M(G;x,y))\big|_{x=1}$ \\
\hline
\end{tabular}
\end{center}
\caption{Three more topological indices and the corresponding formulas}
\label{table2}
\end{table}

%%%%%%%%%%%%%%%%%%%%%%%%%%%%%%%%%%%%%%%%%%%%%%%%%%%%%%%%%%%%%%%%%%%%%%%%%%%%%
\section{How to compute the $M$-polynomial}
%%%%%%%%%%%%%%%%%%%%%%%%%%%%%%%%%%%%%%%%%%%%%%%%%%%%%%%%%%%%%%%%%%%%%%%%%%%%%

As we have demonstrated in the previous section, the computation of the degree based topological indices can be reduced to the computation of the $M$-polynomial. In this section we first give some general remarks how to compute this polynomial and then derive the polynomial for some typical classes of graphs from the literature. In the general setting we follow an approach of Gutman~\cite{gutman-2002} for which some definitions are needed first. 

As we are primarily interested in chemical applications, we will restrict us here to {\em chemical  graphs} which are the connected graphs with maximum degree at most 4. (However, the statements can be generalized to arbitrary graphs.) For a chemical graph $G=(V,E)$, let $n=|V(G)|$, $m=|E(G)|$, and let $n_i$, $1\le i\le 4$, be the number of vertices of degree $i$. Note first that $m_{11} = 0$ whenever $G$ has at least three vertices (and is connected). For the other coefficients $m_{ij}$ of the $M$-polynomial Gutman observed that the following relations hold: 
\begin{eqnarray}
n_1 + n_2 + n_3 + n_4 & = & n \label{eq-3} \\
m_{12} + m_{13} + m_{14} & = & n_1 \label{eq-4} \\
m_{12} + 2m_{22} + m_{23} + m_{24} & = & 2n_2 \label{eq-5} \\
m_{13} + m_{23} + 2m_{33} + m_{34} & = & 3n_3 \label{eq-6} \\
m_{14} + m_{24} +m_{34} + 2m_{44} & = & 4n_4  \label{eq-7}\\
n_1 + 2n_2 + 3n_3 + 4n_4 & = & 2m\,. \label{eq-8}
\end{eqnarray}
These equalities are linearly independent. In some examples, all the quantities $m_{ij}$ (and consequently the $M$-polynomial) can be found directly. On the other hand, one can determine some of the $m_{ij}$'s first and then the remaining ones can be obtained from the above relations. We now demonstrate this on several examples. 

%%%%%%%%%%%%%%%%%%%%%%%%%%%%%%%%%%%%%%%%%%%%%%%%%%%%%%%%%%%%%%%%%%%%%%%%%%%%%
\subsection{Polyomino chains}
%%%%%%%%%%%%%%%%%%%%%%%%%%%%%%%%%%%%%%%%%%%%%%%%%%%%%%%%%%%%%%%%%%%%%%%%%%%%%

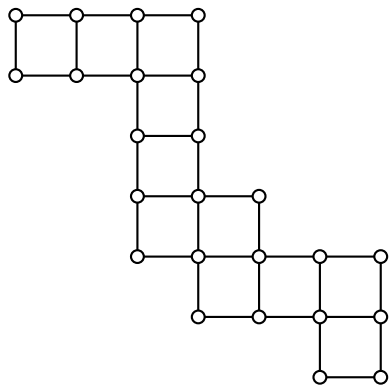
\begin{figure}[ht!]
\begin{center}
\begin{tikzpicture}[scale=0.8,style=thick]
\def\vr{3pt} % \vr = vertex radius;  Set \vr = 2/scale for uniform sizing of vertices
%% edges %%
\draw (0,6) -- (3,6); 
\draw (0,5) -- (3,5); 
\draw (2,4) -- (3,4); 
\draw (2,3) -- (4,3); 
\draw (2,2) -- (6,2); 
\draw (3,1) -- (6,1); 
\draw (5,0) -- (6,0);
\draw (0,6) -- (0,5); 
\draw (1,6) -- (1,5); 
\draw (2,6) -- (2,2); 
\draw (3,6) -- (3,1); 
\draw (4,3) -- (4,1); 
\draw (5,2) -- (5,0); 
\draw (6,2) -- (6,0);
%
%% vertices %%%
\draw (0,6)  [fill=white] circle (\vr);
\draw (1,6)  [fill=white] circle (\vr);
\draw (2,6)  [fill=white] circle (\vr);
\draw (3,6)  [fill=white] circle (\vr);
\draw (0,5)  [fill=white] circle (\vr);
\draw (1,5)  [fill=white] circle (\vr);
\draw (2,5)  [fill=white] circle (\vr);
\draw (3,5)  [fill=white] circle (\vr);
\draw (2,4)  [fill=white] circle (\vr);
\draw (3,4)  [fill=white] circle (\vr);
\draw (2,3)  [fill=white] circle (\vr);
\draw (3,3)  [fill=white] circle (\vr);
\draw (4,3)  [fill=white] circle (\vr);
\draw (2,2)  [fill=white] circle (\vr);
\draw (3,2)  [fill=white] circle (\vr);
\draw (4,2)  [fill=white] circle (\vr);
\draw (5,2)  [fill=white] circle (\vr);
\draw (6,2)  [fill=white] circle (\vr);
\draw (3,1)  [fill=white] circle (\vr);
\draw (4,1)  [fill=white] circle (\vr);
\draw (5,1)  [fill=white] circle (\vr);
\draw (6,1)  [fill=white] circle (\vr);
\draw (5,0)  [fill=white] circle (\vr);
\draw (6,0)  [fill=white] circle (\vr);
%
%% text %%%%%%%%%%%%%%%%%%%%%%%%%%%%%%%%%
% \draw (1.5,0) node {$\cdots$};

\end{tikzpicture}
\end{center}
\caption{Polyomino chain $B(11;6;3,4,2,2,3,2)$}
\label{fig:chain}
\end{figure}

Let $B(n;r;\ell_1,\ldots,\ell_r)$ be a polyomino chain with $n$ squares, arranged in $r$ segments of lengths $\ell_1,\ldots,\ell_r$, $\ell_i\ge 2$, see Fig.~\ref{fig:chain} for  the particular case $B(11;6;3,4,2,2,3,2)$. 

%%%%%%%%%%%%%%%%%%%%%%%%%%%%%%%%%%%%%%%%%%%%%%%%%%%%%%%%%%%%%
%%%%%%%%%%%%%%%%%%%%%% polyomino chain %%%%%%%%%%%%%%%%%%%%%%
%%%%%%%%%%%%%%%%%%%%%%%%%%%%%%%%%%%%%%%%%%%%%%%%%%%%%%%%%%%%%

Let $a$ denote the number of segments of length 2 that occur at the extremities of the chain ($a = 0$, 1, or 2) and  let $b$ denote the number of non-extreme segments of length 2. In Fig.~\ref{fig:chain} we have $a = 1$, $b = 2$. As we shall see, the knowledge of the parameters $n$, $r$, $a$, and $b$ of a chain is sufficient to determine its $M$-polynomial. It is easy to see that $|V| = 2n + 2$ and $|E| = 3n + 1$. The numbers $n_i$ of vertices of degree $i$ are given by  $n_2 = r + 3$ (number of outer corners), $n_4 = r - 1$ (number of inner corners), and $n_3 = |V| - n_2 - n_4 = 2(n - r)$. We have $m_{22} = 2$, $m_{44} = b$ (the edges that halve the non-extreme segments of length 2), and $m_{24} = a + 2b$ (each extreme segment of  length $2$ contributes $1$ and each non-extreme segment of  length $2$ contributes $2$ $24$-edges). Now, taking into account that there are no vertices of degree $1$, from~\eqref{eq-5}-\eqref{eq-7} we obtain the expressions for $m_{23}$, $m_{33}$, and $m_{34}$. Setting $B_n = B(n;r;\ell_1,\ldots,\ell_r)$ this  yields
$$M(B_n; x,y) = 2x^2 y^2 + (2r - a - 2b + 2)x^2 y^3 + (4r - a - 4b - 4)x^3y^4 + bx^4 y^4\,.$$
For example, for the first Zagreb index (corresponding to $f(x,y) = x + y$), either from~\eqref{eq:degree-index-2} or from the first entry of Table~\ref{table1} we obtain $18n +2r - 4$, in agreement with~\cite[Theorem 2.1]{yaas-2012}. Similarly, for the second Zagreb index (corresponding to $f(x,y) = xy$) we obtain $27n + 6r - 19 - a - b$, in agreement with~\cite[Theorem 2.4]{yaas-2012}.  
For the special case of the linear chain $L_n$ $(n\ge 3)$, corresponding to $r = 1$, $a = b = 0$, we obtain
$M(L_n;x,y) = 2x^2 y^2 + 4x^2 y^3 + (3n - 5)x^3 y^3.$
For the special case of the zig-zag chain $Z_n$, corresponding to $r = n-1$, $a = 2$, $b = r - 2 = n -3$, we obtain
$M(Z_n; x,y) = 2x^2 y^2 + 4x^2 y^3 + 2(n - 2)x^2 y^4 + 2x^3 y^4 + (n - 3)x^4 y^4.$
We remark that for the first (second) Zagreb index we obtain $20n - 6$ (resp. $32n - 24$), in disagreement with the slightly mistaken expressions given in~\cite[Corollary 2.2(ii)]{yaas-2012} and ~\cite[Corollary 2.5 (ii)]{yaas-2012}.

%%%%%%%%%%%%%%%%%%%%%%%%%%%%%%%%%%%%%%%%%%%%%%%%%%%%%%%%%%%%%%%%%%%%%%%%%%%%%
\subsection{Starlike trees}
%%%%%%%%%%%%%%%%%%%%%%%%%%%%%%%%%%%%%%%%%%%%%%%%%%%%%%%%%%%%%%%%%%%%%%%%%%%%%

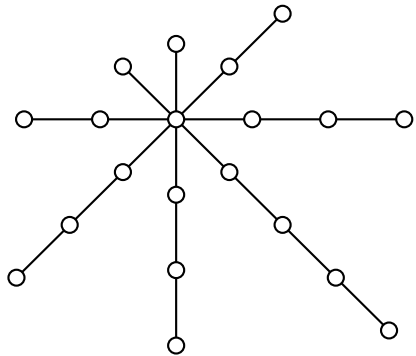
\begin{figure}[ht!]
\begin{center}
\begin{tikzpicture}[scale=1,style=thick]
\def\vr{3pt} % \vr = vertex radius;  Set \vr = 2/scale for uniform sizing of vertices
%% edges %%
\draw (-2,0) -- (3,0); 
\draw (-2.1,-2.1) -- (1.4,1.4); 
\draw (-0.7,0.7) -- (2.8,-2.8); 
\draw (0,1) -- (0,-3); 
%
%% vertices %%%
\draw (0,0)  [fill=white] circle (\vr);
\draw (-1,0)  [fill=white] circle (\vr);
\draw (-2,0)  [fill=white] circle (\vr);
\draw (1,0)  [fill=white] circle (\vr);
\draw (2,0)  [fill=white] circle (\vr);
\draw (3,0)  [fill=white] circle (\vr);
\draw (-0.7,0.7)  [fill=white] circle (\vr);
\draw (0,1)  [fill=white] circle (\vr);
\draw (0.7,0.7)  [fill=white] circle (\vr);
\draw (1.4,1.4)  [fill=white] circle (\vr);
\draw (-0.7,-0.7)  [fill=white] circle (\vr);
\draw (-1.4,-1.4)  [fill=white] circle (\vr);
\draw (-2.1,-2.1)  [fill=white] circle (\vr);
\draw (0,-1)  [fill=white] circle (\vr);
\draw (0,-2)  [fill=white] circle (\vr);
\draw (0,-3)  [fill=white] circle (\vr);
\draw (0.7,-0.7)  [fill=white] circle (\vr);
\draw (1.4,-1.4)  [fill=white] circle (\vr);
\draw (2.1,-2.1)  [fill=white] circle (\vr);
\draw (2.8,-2.8)  [fill=white] circle (\vr);
%
%% text %%%%%%%%%%%%%%%%%%%%%%%%%%%%%%%%%
% \draw (1.5,0) node {$\cdots$};

\end{tikzpicture}
\end{center}
\caption{Starlike tree $S(1,2,3,4,3,3,2,1)$}
\label{fig:star}
\end{figure}

We next consider starlike trees, i.e. trees having exactly one vertex of degree greater than two~\cite{azir-2013}. If the degree of  this vertex (called center) is $n$, then we denote a starlike tree by $S(k_1, \ldots, k_n)$, where the $k_j$'s are positive integers representing the number of edges of the rays emanating from the center, see Fig.~\ref{fig:star} for an example. 

%%%%%%%%%%%%%%%%%%%%%%%%%%%%%%%%%%%%%%%%%%%%%%%%%%%%%%%%%%%%%
%%%%%%%%%%%%%%%%%%%%%%  starlike tree  %%%%%%%%%%%%%%%%%%%%%%
%%%%%%%%%%%%%%%%%%%%%%%%%%%%%%%%%%%%%%%%%%%%%%%%%%%%%%%%%%%%%

Denote $K = \sum_{j=1}^n k_j$ and let $a$ be the number of rays having exactly one edge. The possible vertex degrees are 1, 2, and $n$ and it is easy to see that $m_{11} = 0$, $m_{12} = n - a$, $m_{1,n} = a$, $m_{2,n} = n - a$, $m_{n,n} = 0$. Since the total number of edges is $K$, we have $m_{2,2} = K - (n - a) - a - (n - a) = K + a - 2n$, leading to
$$M(S(k_1, \ldots, k_n); x,y) = (n - a)xy^2 + axy^n + (K + a - 2n)x^2y^2 + (n - a)x^2y^n\,.$$
For example, for the first Zagreb index we obtain (using~\eqref{eq:degree-index-2} or Table~\ref{table1}) $n^2 - 3n + 4K$. Similarly, for the second Zagreb index we obtain $2n^2 - 6n - an + 2a + 4K$. Both results are in agreement with those in~\cite[Corollary 3.7]{azir-2013} obtained in a more complicated manner.

%%%%%%%%%%%%%%%%%%%%%%%%%%%%%%%%%%%%%%%%%%%%%%%%%%%%%%%%%%%%%%%%%%%%%%%%%%%%%
\subsection{Triangulenes}
%%%%%%%%%%%%%%%%%%%%%%%%%%%%%%%%%%%%%%%%%%%%%%%%%%%%%%%%%%%%%%%%%%%%%%%%%%%%% 

\begin{figure}[ht!]
\begin{center}
\begin{tikzpicture}[scale=0.45,style=thick]
\def\vr{4pt} % \vr = vertex radius;  Set \vr = 2/scale for uniform sizing of vertices
% G_1
\draw (-3,0) -- (-2,2) -- (-4,2) -- (-3,0);
\draw (-3,0)  [fill=white] circle (\vr);
\draw (-2,2)  [fill=white] circle (\vr);
\draw (-4,2)  [fill=white] circle (\vr);
\draw (-2.3,0) node {$y_1$};
{\Large \draw (-3,-1) node {$G_1$};}
% G2
\draw (2,0) -- (3,2) -- (1,2) -- (2,0);
\draw (1,2) -- (1.7,4) -- (0.3,4) -- (1,2);
\draw (3,2) -- (3.7,4) -- (2.3,4) -- (3,2);
\draw (2,0)  [fill=white] circle (\vr);
\draw (3,2)  [fill=white] circle (\vr);
\draw (1,2)  [fill=white] circle (\vr);
\draw (1.7,4)  [fill=white] circle (\vr);
\draw (0.3,4)  [fill=white] circle (\vr);
\draw (3.7,4)  [fill=white] circle (\vr);
\draw (2.3,4)  [fill=white] circle (\vr);
\draw (2.7,0) node {$y_2$};
{\Large \draw (2,-1) node {$G_2$};}
% Gn
\draw (9,0) -- (10,2) -- (8,2) -- (9,0);
\draw (8,2) to[out=180,in=270] (6,4)  to[out=90,in=180] (7.5,5.5) to[out=0,in=85] (8.8,3.5) to[out=85,in=0] (8,2);
\draw (10,2) to[out=0,in=270] (12,4)  to[out=90,in=0] (10.5,5.5) to[out=180,in=95](9.2,3.5) to[out=95,in=180] (10,2);
\draw (9,0)  [fill=white] circle (\vr);
\draw (10,2)  [fill=white] circle (\vr);
\draw (10,2)  [fill=white] circle (\vr);
{\Large \draw (9,-1) node {$G_n$};}
\draw (7.3,4.4) node {$G_{n-1}$};
\draw (10.7,4.5) node {$G_{n-1}$};
\draw (7.7,2.6) node {$y_{n-1}$};
\draw (10.4,2.6) node {$y_{n-1}$};
% Tn
\draw (15,3) -- (16,5) -- (17,3) -- (15,3);
\draw (15,3) to[out=180,in=90] (13,1)  to[out=270,in=180] (14.5,-0.5) 
to[out=0,in=275](15.8,1.5) to[out=275,in=0] (15,3);
\draw (17,3) to[out=0,in=90] (19,1)  to[out=270,in=0] (17.5,-0.5) 
to[out=180,in=265] (16.2,1.5) to[out=265,in=180] (17,3);
\draw (16,5) to[out=180,in=270] (14.5,7)  to[out=90,in=180] (16,8.5) 
to[out=0,in=90] (17.5,7) to[out=270,in=0] (16,5);
\draw (15,3)  [fill=white] circle (\vr);
\draw (16,5)  [fill=white] circle (\vr);
\draw (17,3)  [fill=white] circle (\vr);
\draw (16,7.5) node {$G_n$};
\draw (14.5,0.5) node {$G_n$};
\draw (17.5,0.5) node {$G_n$};
\draw (16,5.5) node {$y_{n}$};
\draw (15,2.5) node {$y_{n}$};
\draw (17.2,2.5) node {$y_{n}$};
{\Large \draw (16,-1) node {$T_n$};}
\end{tikzpicture}
\end{center}
\vspace*{-5mm}
\caption{Graphs $G_1$, $G_2$, $G_n$, and $T_n$}
\label{fig:G1-G2-Gn-Tn}
\end{figure}
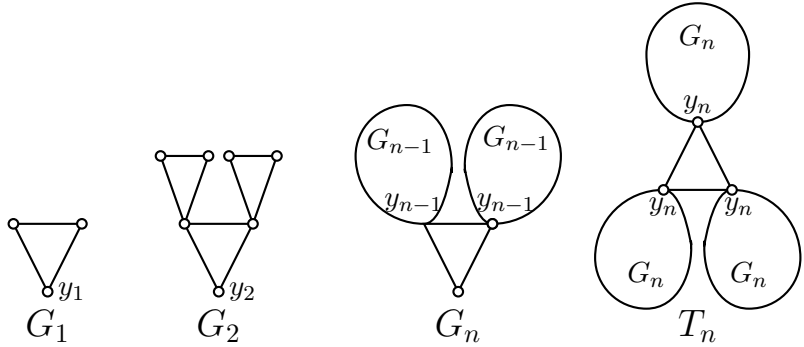

We consider the family of triangulanes $T_n$, defined in Fig.~\ref{fig:G1-G2-Gn-Tn}, where the auxiliary triangulane $G_n$ is defined recursively in the following manner (see also~\cite{as-2011,deutsch-2013,khyo-2008}). Let $G_1$ be a triangle and denote one of its vertices by $y_1$. We define $G_n$ as the circuit of the graphs $G_{n-1}$, $G_{n-1}$, and $K_1$ (the 1-vertex graph) and denote by $y_n$ the vertex where $K_1$ has been placed (see Fig.~\ref{fig:G1-G2-Gn-Tn} again).

%%%%%%%%%%%%%%%%%%%%%%%%%%%%%%%%%%%%%%%%%%%%%%%%%%%%%%%%%%%%%
%%%%%%%%%%%%%%%%%%%%%% G1, G2, Gn, Tn %%%%%%%%%%%%%%%%%%%%%%%%%%
%%%%%%%%%%%%%%%%%%%%%%%%%%%%%%%%%%%%%%%%%%%%%%%%%%%%%%%%%%%%%

The possible vertex degrees of $G_n$ are $2$, $4$ and, using the notation $\mu_{ij}$ instead of $m_{ij}$, it is easy to see that $\mu_{22} = 2^{n-1}$ and $\mu_{24} = 2 + 2^n$. We have $\mu_{22} + \mu_{24} + \mu_{44} =  |E(G_n)| = 3(2^n - 1)$ since $G_n$ consists of $2^n - 1$ triangles. From here we obtain $\mu_{44} = 3\cdot 2^{n-1} - 5$. Thus, as a by-product,
$$M(G_n; x,y) = 2^{n-1} x^2 y^2 + (2 + 2^n) x^2 y^4 + (3\cdot 2^{n-1} - 5) x^4 y^4\,.$$
Now, in order to find the $m_{ij}$'s of $T_n$ from its defining picture, we note that at the attachment of each $G_n$ to the starting triangle, two $24$-edges of $G_n$ turn into two $44$-edges of $T_n$. Consequently,
$$m_{22} = 3\mu_{22} = 3\cdot 2^{n-1},\ m_{24} =  3(\mu_{24} - 2) = 3\cdot 2^n,\ m_{44} = 3(\mu_{44} + 2) + 3 = 3(3\cdot 2^{n-1} - 2)$$
%
%m_{44} = 3(3*2^{n-1} - 2) --->  m_{44} = 3(\mu_{44} + 2) + 3 = 3(3*2^{n-1} - 2)
and
$$M(T_n; x,y) =  3\cdot 2^{n-1} x^2 y^2 + 3\cdot 2^n x^2 y^4 + 3(3\cdot2^{n-1} - 2)x^4 y^4\,.$$

\section*{Acknowledgments}

Work supported in part by the Research Grant P1-0297 of the Ministry of Higher Education, Science and Technology Slovenia.

\end{document}